\documentclass[11pt,a4paper]{article}

\usepackage{amsmath,amssymb}
\usepackage{graphicx}
\usepackage{multirow}
\usepackage{url}

\newcommand{\RR}{\mathbb{R}}
\newcommand{\norm}[1]{\|#1\|}
\DeclareMathOperator{\OT}{OT}

\begin{document}

\title{DOTmark --- A Benchmark for Discrete Optimal Transport}

\author{J\"orn Schrieber\thanks{Supported by the DFG Research Training Group 2088 ``Discovering Structure in Complex Data: Statistics Meets Optimization and Inverse Problems''.}\ \thanks{Institute for Mathematical Stochastics, University of G\"ottingen, Goldschmidtstr. 7, 37077 G\"ottingen, Germany. Email addresses: joern.schrieber-1@mathematik.uni-goettingen.de, dominic.schuhmacher@mathematik.uni-goettingen.de and gottschlich@math.uni-goettingen.de}\, , \ 
	Dominic Schuhmacher\footnotemark[2] \ and
        Carsten Gottschlich\footnotemark[2] \\[2.5mm]
        University of G\"ottingen
}

\maketitle

\begin{abstract}
The Wasserstein metric or earth mover's distance (EMD) is a useful tool in statistics, machine learning and computer science with many applications to biological or medical imaging, among others. Especially in the light of increasingly complex data, the computation of these distances via optimal transport is often the limiting factor. Inspired by this challenge, a variety of new approaches to optimal transport has been proposed in recent years and along with these new methods comes the need for a meaningful comparison.

In this paper, we introduce a benchmark for discrete optimal transport, called DOTmark, which is designed to serve as a neutral collection of problems, where discrete optimal transport methods can be tested, compared to one another, and brought to their limits on large-scale instances. It consists of a variety of grayscale images, in various resolutions and classes, such as several types of randomly generated images, classical test images and real data from microscopy.

Along with the DOTmark we present a survey and a performance test for a cross section of established methods ranging from more traditional algorithms, such as the transportation simplex, to recently developed approaches, such as the shielding neighborhood method, and including also a comparison with commercial solvers. \\[2mm]
\end{abstract}

\section{Introduction}
\label{intro}

Despite being a classical problem, optimal transport appears in a plethora of modern applications, such as
image retrieval \cite{RubnerEtAl2000}, 
phishing web page detection~\cite{FuWenyinDeng2006}, 
measuring plant color differences \cite{KendalEtal2013} 
or shape matching \cite{GraumanDarrell2004}. 
Motivated by the availability of ever larger and more complex data, which often can be cast into the form of one or several measures, and the corresponding development of analytical methods for such data, came the need for more efficient ways of computing optimal transportation plans and evaluating their costs. Accordingly, the last few years have seen the advent of many new methods for obtaining or approximating solutions to large transport problems. 

While many of these methods undoubtedly mean substantial progress compared to what was available a decade ago, it is nearly impossible from the current literature to figure out how various of these methods compare to one another and which method is most suitable for a given task. This is mainly due to the fact that only for a few of the modern methods user-friendly code is publicly available. What is more, many articles that introduce new methods compare their computational performance only on a restricted set of self-generated ad-hoc examples and typically demonstrate improved performance only in comparison to some classical method or to a simplified version that does not use the novelty introduced.

The purpose of the present article is twofold. First, we propose a collection of real and simulated images, the DOTmark, that is designed to span a wide range of different mass distributions and serves as a benchmark for testing optimal transport algorithms. The data can be downloaded at \url{www.stochastik.math.uni-goettingen.de/DOTmark/}. We invite other researchers to use this benchmark, report their results, and thus help building a more transparent picture of the suitability of different methods for various tasks.  

The second purpose is to provide a survey and a performance test based on the \mbox{DOTmark} for a cross section of established methods. Since not much code is freely available, we have used previous implementations of our own (done to the best of our knowledge) of various methods, added an implementation of the recently proposed shielding neighborhood method \cite{Schmitzer2016} and let them compete against each other. This also allows us to draw conclusions about the behavior of different methods on different types of input data.
In order to make this comparison as meaningful as possible, we restricted ourselves to using only singlescale methods and the squared Euclidean distance as a cost function. We hope this comparison will provide a first spark for a healthy competition of various methods in the public discussion.

\section{Brief Theoretical Background}
\label{theo}

For the present context it is sufficient to restrict ourselves to optimal transport on $\RR^d$. Let $X,Y$ be subsets of $\RR^d$ and let $\mu$ and $\nu$ be probability measures on $X$ and $Y$, respectively. In this paper we will always have $X=Y$, but using different notation for domain and target space makes definitions easier to grasp.

\begin{figure}[ht]

\begin{minipage}{0.33\textwidth}
  \includegraphics[width=0.975\textwidth]{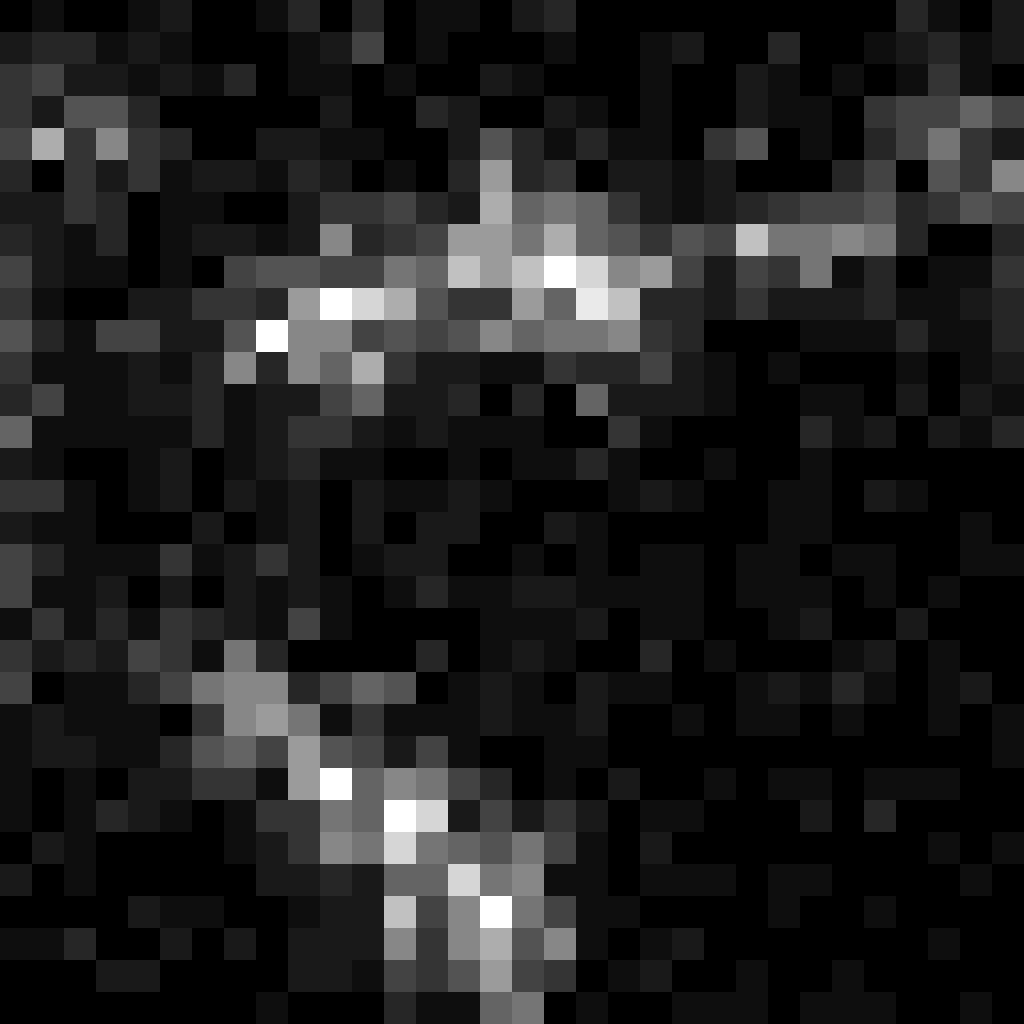} \vspace{6pt}\\
  \includegraphics[width=0.975\textwidth]{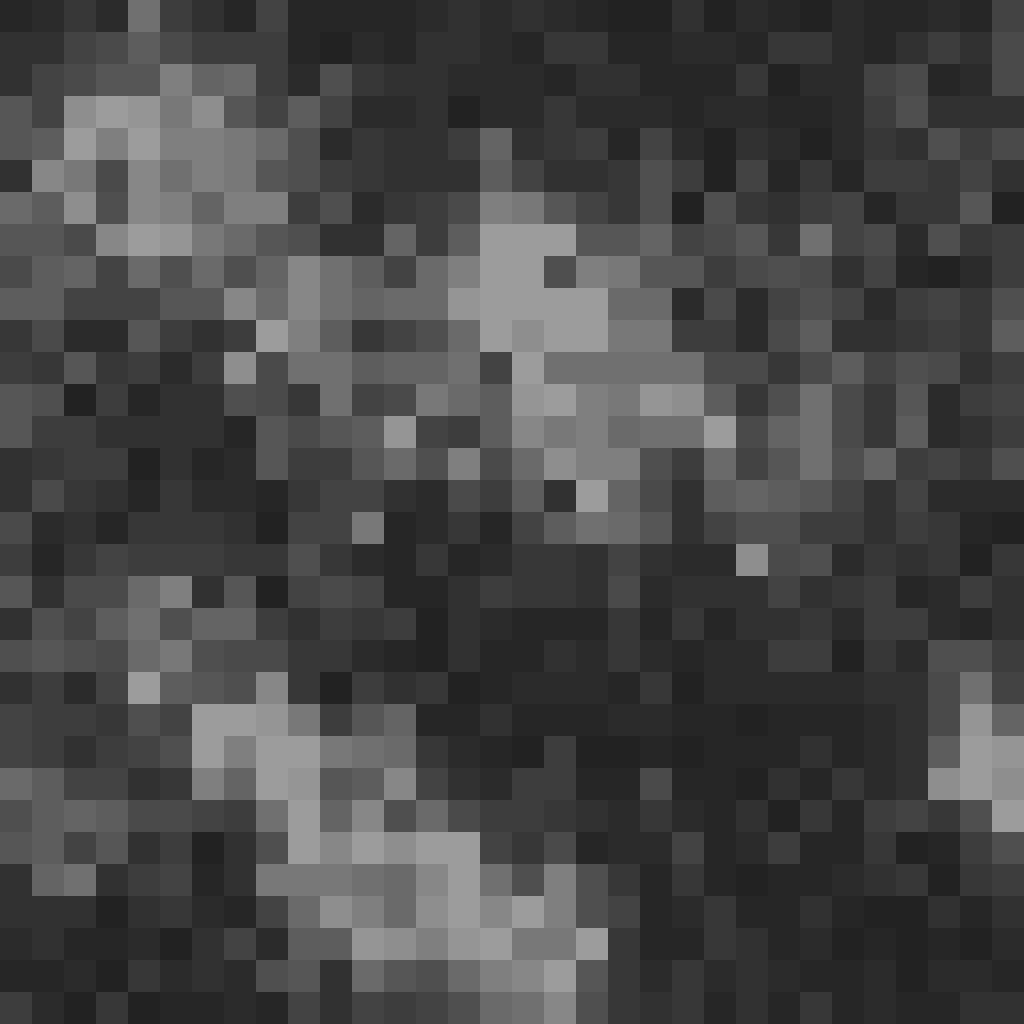}
\end{minipage}
\begin{minipage}{0.66\textwidth}
 \includegraphics[width=\textwidth]{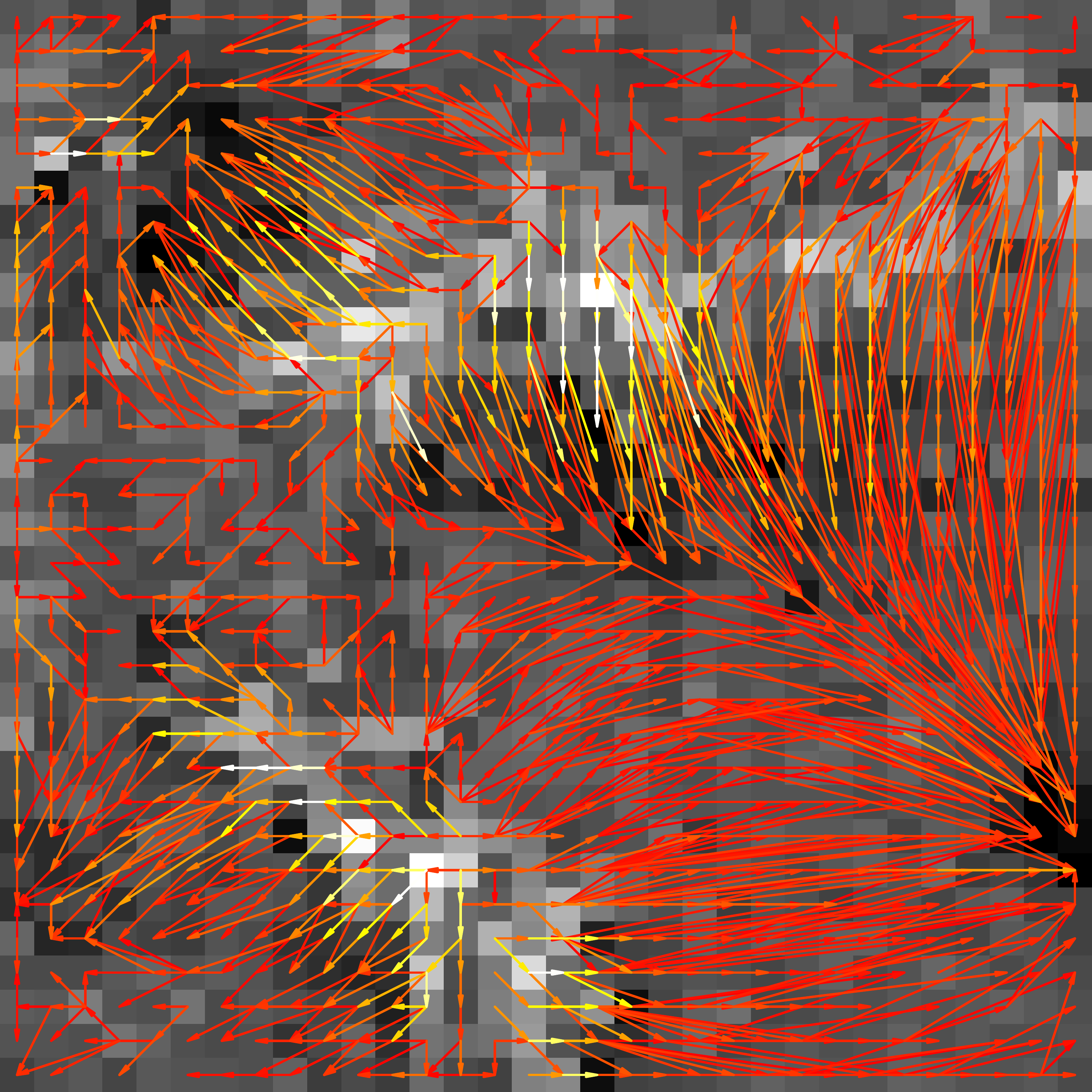}
\end{minipage}

  \caption{\emph{Left panels:} Two tiny clippings $A$ (top) and $B$ (bottom) from STED microscopy images of mitochondrial networks. \emph{Right panel:} The difference $A-B$ of the first two panels with an optimal transference plan for $p=2$ superimposed. Arrows show from where to where mass is transported in the optimal transference plan. The colors indicate the amount of mass from dark red (small) to bright yellow (large). Since mass from individual sites is split (indicated by several arrows leaving the site), this transference plan cannot be represented by a transport map.}
  \label{fig:microdemo}
\end{figure}

A \emph{transport map} $T$ is any (measurable) map $X \to Y$ that transforms the measure $\mu$ into the measure $\nu$. More precisely it satisfies $\mu(T^{-1}(B)) = \nu(B)$ for every measurable $B \subset Y$. A \emph{transference plan} is a measure $\pi$ on $X \times Y$ with marginals $\pi(\cdot \times Y) = \mu$ and $\pi(X \times \cdot) = \nu$. The set of transference plans from $\mu$ to $\nu$ is denoted by $\Pi(\mu,\nu)$. Any transport map $T$ from $\mu$ to $\nu$ defines a transference plan $\pi_T$ from $\mu$ to $\nu$ as the unique measure satisfying $\pi_T(A \times B) = \mu(A \cap T^{-1}(B))$ for all measurable $A \subset X$ and $B \subset Y$. Not every transference plan $\pi$ can be represented in this way, because transference plans allow mass from one site $x \in X$ to be split between multiple destinations, which is not possible under a transport map. Figure \ref{fig:microdemo} shows such an example.

We assume that the cost of transporting a unit mass from $x \in X$ to $y \in Y$ is $c_p(x,y) = \norm{x-y}^p$ for some $p \geq 1$. The minimum cost for transferring $\mu$ to $\nu$ is then given by
\begin{equation}
\label{eq:mincost}
  C_p(\mu,\nu) = \min_{\pi \in \Pi(\mu,\nu)} \int_{X \times Y} \norm{x-y}^p \; d\pi(x,y).
\end{equation}

Taking the $p$-th root, we obtain the \emph{Wasserstein metric} $W_p$. More precisely we have 
\begin{equation*}
  W_p(\mu,\nu) = C_p(\mu,\nu)^{1/p}
\end{equation*}
for any measures $\mu$ and $\nu$ that satisfy $\int_X \norm{x}^p \, d\mu(x) < \infty$ and $\int_Y \norm{y}^p \, d\nu(y) < \infty$. In order to evaluate the Wasserstein metric, we need to find an optimal solution to (\ref{eq:mincost}), i.e., a minimizing transference plan $\pi$. This problem is often referred to as the Kantorovich formulation of optimal transport. Note that by Theorem~4.1 in~\cite{Villani2009} a minimizing $\pi$ always exists. However, it neither has to be unique nor representable in terms of an optimal transport map.

Often one would like to compare data sets that are available as images from a certain source, e.g.\ real photography, astronomical imagery, or microscopy data. We may think of such images as discrete measures on a grid. For example, the first two panels in Figure~\ref{fig:microdemo} show tiny clippings from STED microscopy images of mitochondrial networks. A question of interest might be whether both images stem from the same part of the network, which can in principle be answered by finding an optimal transference plan (third panel in Figure~\ref{fig:microdemo}) and computing the Wasserstein distance. Note that this coarse resolution is not representative for a serious analysis, but was only chosen for illustrative purposes.

Even if the measures we would like to consider are more general probability measures, we can always approximate them (weakly) by a discrete probability measure, e.g.\ by considering the empirical distribution of a sample from the general measure or based on a more sophisticated quantization scheme. Lemma~8.3 in \cite{BickelFreedman1981} characterizes when optimal costs are approximated in this way (e.g.\ always if $X$ and $Y$ are compact). Theorem~5.20 in~\cite{Villani2009} and the subsequent discussion give sufficient conditions about the approximation of optimal transference plans. 

Assume now that we have discrete measures of the form $\mu = \sum_{i=1}^{m} \mu_i \delta_{x_i}$ and $\nu = \sum_{j=1}^{n} \nu_j \delta_{y_j}$ and write $c_{ij} = \norm{x_i-y_j}^p$. In what follows, we always have $m = n$, and $(x_i)_{1 \leq i \leq m} = (y_j)_{1 \leq j \leq n}$ form a regular square grid in $\RR^2$, but since it is more intuitive, we keep different notation for source locations and target locations. Let $\pi_{ij}$ be the amount of mass transported from $x_i$ to $y_j$. Then, the problem \eqref{eq:mincost} can be rewritten as a linear program:
\begin{alignat*}{3}
	\OT \quad &&\min \quad &\sum_{i = 1}^m \sum_{j = 1}^n c_{ij} \pi_{ij}\\
 	&&\text{subject to} \quad & \sum_{j=1}^n \pi_{ij} = \mu_i \ &&\forall i = 1, \dots, m \\
 	&& &\sum_{i=1}^m \pi_{ij} = \nu_j \ &&\forall j = 1, \dots, n\\
 	&& & \pi_{ij} \ge 0
\end{alignat*}
This is the classic transportation problem from linear programming. Efficient ways of solving this problem for small to medium sized ($m$ and) $n$ have been known since the middle of the last century. However, in the context of modern optimal transport problems it has become necessary to solve such problems efficiently at a scale where ($m$ and) $n$ are many thousands or even tens of thousands and more. Currently, this cannot be done with the classical algorithms and requires to utilize the geometry of the problem in one way or the other.

\section{Benchmark}
\label{bench}
Our philosophy in compiling this benchmark was to represent a wide range of theoretically different structures, while incorporating typical images that are used in praxis and/or have been used for previous performance tests in the literature. We refer to it as DOTmark, where DOT stands for discrete optimal transport.

\begin{table}[th]
\centering
\renewcommand{\arraystretch}{1.2}
\begin{tabular}{r|l|p{0.667\textwidth}|}
\textbf{\#} & \textbf{Name} & \textbf{Description} \\
\hline \hline
1 & WhiteNoise & i.i.d.\ uniformly distributed values in $[0,1]$ at each pixel \\ \hline
2 & GRFrough & GRF with $\sigma^2=1$, $\nu = 0.25$, $\gamma = 0.05$ \\ \hline
3 & GRFmoderate & GRF with $\sigma^2=1$, $\nu = 1$, $\gamma = 0.15$ \\ \hline
4 & GRFsmooth & GRF with $\sigma^2=1$, $\nu = 2.5$, $\gamma = 0.3$ \\ \hline
5 & LogGRF & $\exp$-function of a GRF with $\sigma^2=1$, $\nu = 0.5$, $ \gamma = 0.4$ \\ \hline
6 & LogitGRF & Logistic function of a GRF with $\sigma^2=4$, $\nu = 4.5$, $\gamma = 0.1$ \\ \hline
7 & CauchyDensity & Bivariate Cauchy density with random center and a varying scale ellipse \\ \hline
8 & Shapes & An ad-hoc choice of simple geometric shapes \\ \hline
9 & ClassicImages & Standard grayscale test images used in image processing \\ \hline
10 & Microscopy & Clippings from STED microscopy images of mitochondria \\
\end{tabular}
\caption{The 10 classes in the DOTmark with details about their creation. GRF stands for Gaussian random field. For technical details and the meaning of the parameters see text.}
\label{tab:10classes}
\end{table}

The benchmark consists of 10 classes of 10 different images (in what follows sometimes called mass distributions or measures), each of which is available at the 5 different resolutions from $32 \times 32$ to $512 \times 512$ (in doubling steps per dimension). This allows for a total of 45 computations of Wasserstein distances between two images for any one class at any fixed resolution. Table~\ref{tab:10classes} gives an overview of how the classes were created. Classes 1--7 are random simulations of scenarios based on various probability distributions. Images at different resolutions are generated independently from each other but according to the same laws. Classes 8--10 were obtained by ad-hoc choices of simple geometric shapes, classic test images and 
images of mitochondria acquired using STED super-resolution microscopy \cite{IllgenEtAl2014,JansEtAl2013,WurmEtAl2011}.
For geometric shapes and classic test images the various resolutions available are coarsenings of a single image. 
For the microscopy images different clippings of various sizes have been selected from larger images to obtain the various resolutions.

We shifted and scaled the pixel values for all classes and randomly redistributed a small percentage of the mass in order to achieve non-negative integer values at each pixel with an average of $10^5$. We chose integer values to make the benchmark (directly) accessible to a wide range of algorithms and to be able to verify correctness of the optimal transport cost precisely, at least in the case $p=2$, where integer transportation costs between grid points may be assumed.

\begin{figure}[th]
  \hspace*{-2.5pt}\includegraphics[width=1.006452\textwidth]{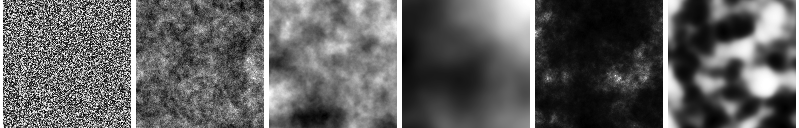}\\[1mm]
  \hspace*{-2.5pt}\includegraphics[width=1.006452\textwidth]{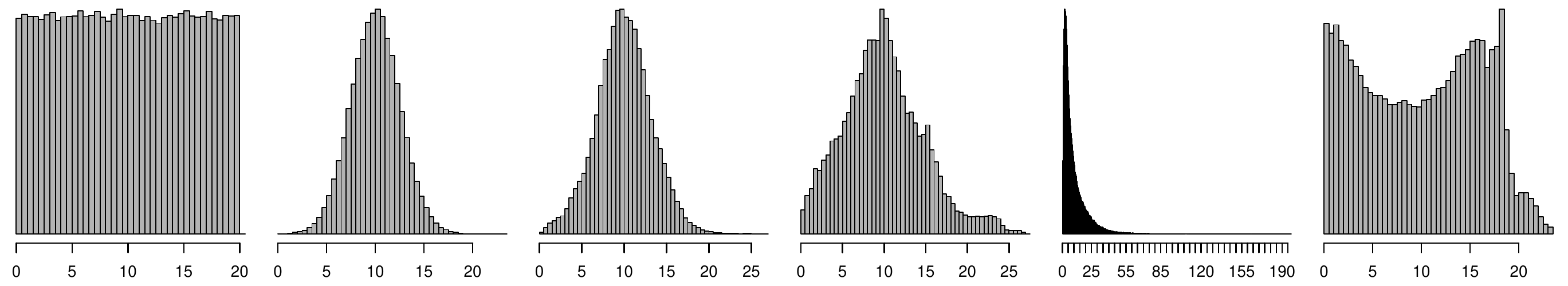}
  \vspace*{-6mm}

  \caption{\emph{Top row:} One representative at resolution $128 \times 128$ for each of the completely randomly generated classes 1--6.
  \emph{Bottom row:} Average histograms of all images at resolution $128 \times 128$ in classes 1--6. A regular bin width of 5000 was used. The annotated pixel values are in multiples of $10^4$.}
  \label{fig:bmrandom}
\end{figure}

\begin{figure}[th]
  \includegraphics[width=\textwidth]{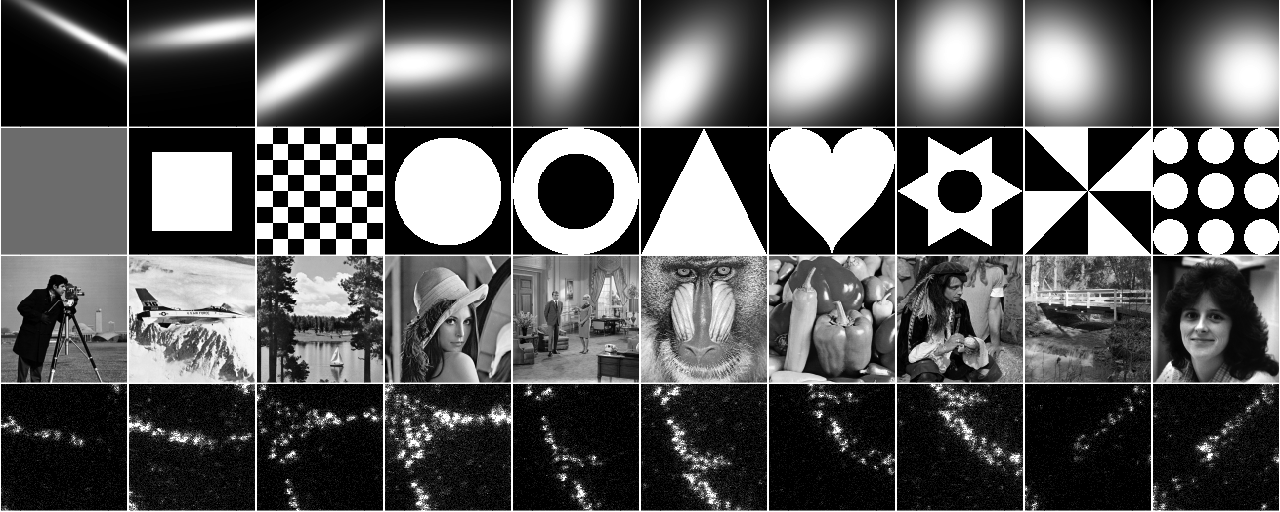}
  \caption{The images in the classes 7--10 at resolution $128 \times 128$.}
  \label{fig:bmfixed}    
\end{figure}

Figure~\ref{fig:bmrandom} shows the first image of each of the classes 1--6 along with average histogram over all members of the class. Figure~\ref{fig:bmfixed} shows the complete collection of images in classes 7--10.

We provide some more details on how classes 2--6 are generated. The images are simulated from stationary centered Gaussian random fields (GRF) on $[0,1]^2$ with Mat\'ern covariance function $k := k_{\sigma^2,\nu,\gamma} : \RR^2 \times \RR^2 \to \RR$,
\begin{equation*}
  k_{\sigma^2,\nu,\gamma}(x,y) = \sigma^2 \frac{2^{\nu-1}}{\Gamma(\nu)} \biggl( \sqrt{2\nu} \frac{\norm{x-y}}{\gamma} \biggr)^{\nu} K_{\nu} \bigl( \sqrt{2\nu} \norm{x-y}/\gamma \bigr),
\end{equation*}
where $K_{\nu}$ is the modified Bessel function of the second kind of order $\nu$. In brief this means that the pixel values are distributed according to a multivariate normal distribution with mean vector zero and covariance matrix $(k(x_i,x_j))_{1 \leq i,j \leq m}$, where $x_i$, $1 \leq i \leq m$, is an enumeration of the pixel centers. The Mat\'ern covariance function is a popular choice in spatial statistics. Its significance comes from the fact that in addition to having parameters $\sigma^2>0$ for the variance and $\gamma>0$ for the range of the covariance, it also has a parameter $\nu>0$ that allows to control the regularity of the random image created from very rough ($\nu$ small) to very smooth ($\nu$ large). Accordingly, classes 2--4 go from very rough with short range dependence to quite smooth with long range dependence. Class~5 is rough with long range dependence, which is hard to see from Figure~\ref{fig:bmrandom} because of the exponential function applied to the pixel values. Class~6 is very smooth with medium range dependence and the logistic function $\psi: \RR \to [0,1]$, $\psi(x) = e^x \bigl/ (1+e^x)$ was applied to the pixel values. See~\cite{GneitingGuttorp2010} for more theoretical details about Gaussian random fields. Our simulations were performed using the {\sf R} package \texttt{RandomFields}; see~\cite{SchlatherEtAl2016} and~\cite{R}.

Histograms 4 and 6 deviate quite a bit from the theoretical histograms expected due to the rescaling and redistribution of mass that we apply in order to obtain mass distributions that are non-negative, integer-valued, and have an average of $10^5$. Also, due to long range dependence and smoothness, histogram 4 is based on a sample of much smaller effective size from the normal distribution than histograms 2 and 3.

On the whole we consider this a reasonable and versatile benchmark for many (planar and typically grid-based) optimal transport algorithms. It covers a wide selection of types of mass distributions whose comparisons are useful for theoretical or practical investigations. Similar types have been considered in the literature before. Gottschlich and Schuhmacher \cite{GottschlichSchuhmacher2014} have considered a sparse version of the WhiteNoise class. Schmitzer \cite{Schmitzer2016} considered sums of randomly scaled and positioned Gaussian densities (sometimes filtered by discontinous masks), which is a somewhat different type of random function generation along the lines of our classes 2--7. Further random or deterministic functions and geometric shapes were considered by Benamou, Carlier, Cuturi, Froese, Nenna, Oberman, Peyr{\'e}, Ruan and others; see e.g.\ \cite{BFO2014}, \cite{BCCNP2015}, \cite{ObermanRuan2015}. The use of real grayscale images as in class~9, but also color images, is abundant in the computer science literature (e.g.\ \cite{RubnerEtAl2000}, \cite{PeleWerman2009}), where optimal transport is typically considered on some feature space. M\'erigot \cite{Merigot2011} illustrated and tested his algorithm on grayscale images directly. Biological imagery has been used in \cite{GottschlichHuckemann2015} (fingerprints in feature space), \cite{GerberMaggioni2015} (brain MRIs) and elsewhere.

Another example that is frequently used, in particular in the statistics and machine learning literature, are images from the MNIST handwritten digit database. Due to the low resolution of these images we do not consider them, but we might include other images of handwritten text in later revisions of the benchmark.

All in all, we see the current benchmark as a first stable version. We are happy to adapt and extend future versions based on feedback from other researchers.

\section{Tested Methods}
\label{meth}

In what follows, we describe the methods that we have tested on the DOTmark. Due to the large number of suggested methods, which unfortunately is not well reflected in the number of user-friendly implementations available, some restrictions had to be made. We chose methods with a good track record, such as the AHA method, as well as some new and promising methods, such as the shielding method. In order to make our comparison as meaningful as possible, we abstained from using methods of a purely approximative nature, such as Sinkhorn scaling~\cite{Cuturi2013}. 

For all of the tested methods there exist multiscale versions. Sometimes these are tailor-made, like for the AHA and the shielding methods (see~\cite{Merigot2011} and~\cite{Schmitzer2016}). Sometimes these are just relatively simple but efficient procedures exploiting the fact that all mass distributions considered live on a square grid in $\RR^2$. Such a simple strategy may be as follows: First solve the transport problem on coarsened images (e.g.\ obtained by adding up the pixel values in contiguous squares of four pixels); then refine the obtained transport plan in a suitable way so that a feasible transport plan for the finer images is obtained; finally solve the original (fine) problem using this transport plan as a starting value.

In our experience this simple strategy already results in an improvement by a factor of 2 to 5 in the transportation simplex at resolution $64 \times 64$. A more elaborate, efficient, but not entirely rigorous alternative was proposed in \cite{ObermanRuan2015}. Since the precise variant and implementation of a multiscale method may distort competition and distract from the merit of an algorithm as such, we decided not to use \emph{any} multiscale methods for this first comparison.

\subsection{Transportation Simplex}
\label{tps}
One of the classical optimal transport methods we test in this paper is the transportation simplex, sometimes also referred to as the revised simplex. It is a specialized version of the network simplex and described in detail for example in \cite{LuenbergerYe2008}. Like other simplex variants, the transportation simplex has two phases: one phase to construct an initial basic feasible solution and another phase to improve this solution to optimality. Typically, the majority of time is spent in the second phase, as an initial solution to optimal transport is easily obtained.

There are quite a few different ways to construct an initial basic feasible solution. For a selection, see \cite{GottschlichSchuhmacher2014}. The method we use here, the modified row minimum rule, has a universally solid performance both in runtime and quality of the solution constructed. We iterate over all source locations (rows) $x_i \in X$ that still have mass left and choose for each source the available transport $\pi_{ij}$ with the least cost and include it with the maximal amount possible in our solution. This process is repeated until all sources are depleted. The solution we obtain is automatically basic.

If we interpret the source and target locations as nodes in a graph and draw arcs for every possible transport, we get a complete bipartite graph. Every (non-degenerate) basic feasible solution can now be represented by a spanning tree in this graph by choosing all the arcs belonging to active transports in our solution. Given a basic feasible solution, a simplex step is performed as follows:
\begin{itemize}
	\item A new variable (transport $\pi_{ij}$) is selected to enter the basis.
	\item This creates a cycle in the previous tree, which is then identified.
	\item The maximal amount of mass possible is shifted along this cycle, i.e.,\ alternately added to and subtracted from consecutive transports.
	\item A variable that has become zero in the process is removed from the basis.
\end{itemize}

A row minimum strategy is used when searching for a new basic variable: We use the current solution to compute the values of dual variables $u_i$ for each source $x_i$ and $v_j$ for each target $y_j$. Then we scan the non-basic variables $\pi_{ij}$ row by row and compute the reduced costs $r_{ij} = c_{ij} - u_i - v_j$. If we encounter a variable with negative reduced costs, we stop at the end of that row and choose the variable with the least reduced cost encountered thus far as a new basic variable. If no candidates are found among all rows, the current solution is optimal and we stop.

The second phase of the transportation simplex is very similar to the network simplex. The only difference is that we always have a complete bipartite network in the optimal transport problem. In phase one, however, this structure is very beneficial and allows for easy construction of a feasible solution, whereas in the more general network case the introduction of artificial variables is often necessary. More details on the network simplex can also be found in \cite{LuenbergerYe2008}.

In the test we use our implementation of the transportation simplex in Java.

\subsection{Shortlist Method}
\label{short}
At its core, the shortlist method is a variant of the transportation simplex. It comes with three parameters:
\begin{itemize}
	\item A parameter $s$ for the shortlist length.
	\item Parameters $p$ and $k$ that control how many variables are searched to find a new basic variable.
\end{itemize}

Before the optimization starts, a shortlist is created for every source, consisting of the $s$ targets with least transport costs, ordered by cost. The basic feasible solution is again constructed by a modified row minimum rule, where the shortlists are prioritized. After that, the solution is improved by simplex steps similar to the transportation simplex, but the search is limited to the shortlists. The lists are scanned, until either $k$ variables with negative reduced costs are found or $p$ percent of the shortlists have been searched. Then the candidate with least reduced costs is chosen to enter the basis. If no improvement can be achieved within the shortlists, the last solution is improved to global optimality by the same simplex steps as in the transportation simplex. For further details, see \cite{GottschlichSchuhmacher2014}.

Just as with the transportation simplex, we use our Java implementation of the shortlist method for the benchmark. Many routines are shared with the transportation simplex. We use the default parameters presented in \cite{GottschlichSchuhmacher2014}. They were chosen with regard to the problems considered in that paper -- a version of the class WhiteNoise with the Euclidean distance as cost function ($p = 1$) and irregular source and target locations.

\subsection{Shielding Neighborhood Method}
\label{shield}
The shielding neighborhood method, or shortcut method, was introduced by Schmitzer in \cite{Schmitzer2016}. Its main idea is to solve a sequence of sparse (i.e.\ restricted) optimal transport instances instead of the dense (full) problem. The algorithm is proposed as a multiscale method, but for the reasons stated above we use the singlescale variant, which basically works in the same way:

Starting with a basic feasible solution generated with the modified row minimum rule a \textit{shielding neighborhood} for that solution is constructed as described in \cite{Schmitzer2016} for the squared Euclidean distance as cost function. This neighborhood is a small subset $N \subseteq X \times Y$ of the product space and imposes a restricted instance of the problem by only considering transport variables $\pi_{ij}$ such that $(x_i, y_j) \in N$. Due to the significant reduction in the amount of variables, this instance is much faster to solve.

The idea behind the shielding neighborhood is the so-called \textit{shielding condition}, which ensures that for $(x_i, y_j) \not \in N$ a shortcut exists, i.e.,\ a sequence of transports through the neighborhood whose combined cost is not higher than the cost of the direct transport. The algorithm alternates between optimizing the sparse instance of the problem and generating a new shielding neighborhood for the current solution. If a solution is optimal for two successive shielding neighborhoods the algorithm stops. In \cite{Schmitzer2016} it is proved that such a solution is always globally optimal.

For the internal sparse instances we follow the recommendation of that paper and use the network simplex solver of CPLEX with a warm start of the previously optimal basis. The Java API is used to create the models and we implemented the shielding neighborhood generation routine in Java as well. One call of the CPLEX solver for a single sparse instance is what we refer to as one iteration of the shielding method in the remainder of this paper.

Schmitzer published his own code of the method on his website\footnote{\url{https://www.ceremade.dauphine.fr/~schmitzer/}}. But since we encountered difficulties in getting it to run on our server and since we already had our own well-working implementation of this method in place, we decided to use the latter in our test. As all the other methods it was implemented to the best of our knowledge. Since the majority of the runtime is occupied by the internal CPLEX solver, we expect our code to have a similar runtime as the code provided by Schmitzer.

Additionally, we tested the shielding method with an adapted implementation of the transportation simplex from Subsection~\ref{tps} as internal solver, using the same pivot strategy and simplex step routines as before.

\subsection{AHA Method}
\label{aha}
The acronym AHA stands for Aurenhammer, Hoffmann and Aronov, who showed in their seminal paper \cite{AHA1998} that the transport problem with squared Euclidean cost is equivalent to an unrestricted continuous minimization problem for a certain convex objective function $\Phi$.

Our concrete implementation is largely based on \cite{Merigot2011}, but as with other algorithms we only use the singlescale version for comparability. The algorithm computes the optimal transport from an absolutely continuous measure $\mu$ on $\RR^2$ to a discrete measure $\nu$ on $\RR^2$, a problem sometimes referred to as \emph{semidiscrete optimal transport}. The key observation utilized by the AHA method is that any power diagram (a.k.a.\ Laguerre tesselation) governed by the support points $y_1,\ldots,y_n$ of $\nu$ and arbitrary weights $w_1, \ldots, w_n \in \RR$ characterizes an optimal transport plan from $\mu$ to \emph{some} measure living on these support points.
By minimizing the function $\Phi$ in the weights $w_1,\ldots,w_n$ we can find a power diagram that defines an optimal transport to \emph{the correct} measure $\nu$.

Evaluating $\Phi$ at a given weight vector $w \in \RR^n$ involves computation of the power diagram for $w$ and a rather simple integration procedure over each power cell. The gradient of $\Phi$ is accessible, and its $i$-th component is in fact just the difference between the $\nu$-mass at the $i$-th support point and the $\mu$-mass transported to this point under the current power diagram. A Hessian of $\Phi$ is not accessible. We follow \cite{Merigot2011} in using the L-BFGS-B algorithm with Mor\'e--Thuente type line search.
Since we use a continuous optimization method, we typically cannot reach a weight vector $w_{*}$ where the gradient of $\Phi$ is exactly zero (and hence the image measure of $\mu$ is exactly $\nu$), but have to stop when its length is still slightly positive. We thus commit a small controllable error, which we refer to as \emph{precision error}.

In order to make the algorithm applicable to the fully discrete situation we study in the other algorithms, we turn the first image into an absolutely continuous measure $\mu$ by interpreting pixel values as masses uniformly distributed over the squared areas represented by the pixels, rather than centered at grid points. Compared to the other methods this leads to slightly different results, a discrepancy we refer to as \emph{blurring error}.\footnote{Note that the term ``error'' is subjective. We might as well declare that we want to solve the semidiscrete problem, in which case all the other methods commit a ``concentration error''.} 

Unlike for the other methods we use an implementation that is written mainly in C with some minimal {\sf R} overhead. The construction of power diagrams was reimplemented in C based on ideas from the CGAL Regular\_triangulation\_2 package and other sources. For the L-BFGS-B algorithm we used the implementation in the {\sf R} function \texttt{optim}.

\subsection{Solvers}
\label{solvers}
For representatives of LP solvers, we used CPLEX and Gurobi. For both we modeled the optimal transport problem as an LP and used the default parameters. This is denoted CPLEX-Def and Gurobi-Def, respectively. Additionally, we tested the network simplex solver CPLEX provides (CPLEX-NWS). Gurobi does apparently not come with a network solver, but as the Gurobi documentation page for methods\footnote{\url{http://www.gurobi.com/documentation/6.5/refman/method.html#parameter:Method}} recommends the dual simplex for memory intensive models, we included it in our tests (Gurobi-DS). All models were set up using the Java APIs of the solvers.

\section{Computational Results}
\label{comp}
All of our tests were performed using a single core on a Linux server (AMD Opteron Processor 6140 from 2011 with 2.6 GHz).  
Note that much better absolute runtimes can be achieved when using modern CPU hardware. For many of the algorithms considerable further improvements are possible by multithreaded implementations that use multiple CPU cores simultaneously.

In our experiments we placed the main emphasis on ensuring that the commercial solvers and all of our own implementations were run under the same conditions. In particular, they were all restricted to use only one of 32 available cores, which was realized by the Linux kernel feature cgroups.

Pairing any two of the $10$ images in each of the $10$ classes gives $45$ transport problems (``instances'') per class, yielding $450$ instances in total. These were all solved at resolutions $32 \times 32$ and $64 \times 64$ by each of the described methods using the squared Euclidean metric as cost funtion. All optimal transport costs returned were checked for correctness. The AHA method is the only procedure, where we cannot expect precisely correct results due to the errors described in Subsection \ref{aha}. These errors are reported in Subsection~\ref{errors}. All other errors were zero. 

\subsection{Runtimes}
\label{runtimes}
\begin{figure}[p]
  \vspace*{-8mm}
  
  \centering
  \hspace*{-4mm}
  \includegraphics{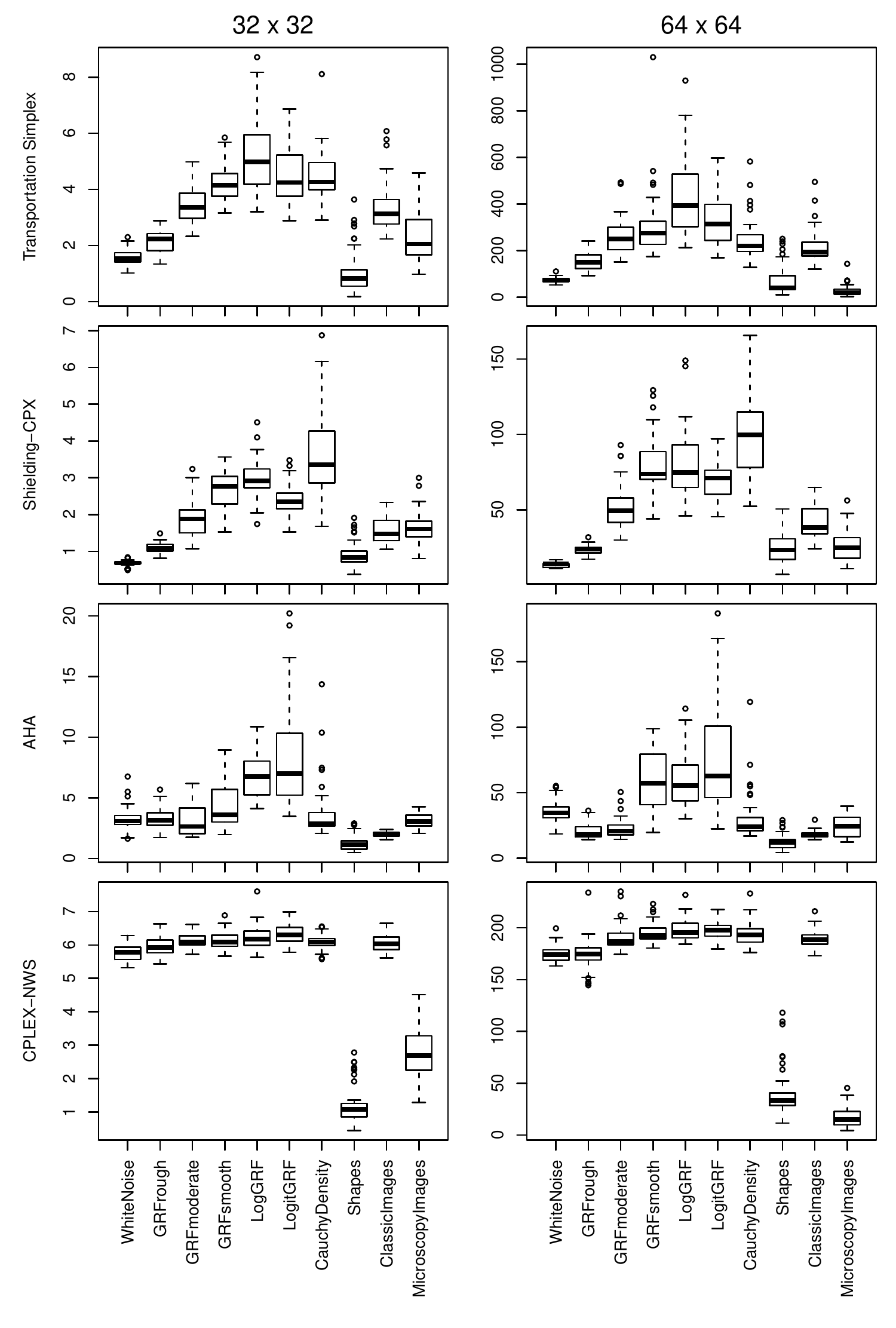}
  \vspace*{-9mm}

  \caption{Boxplots of the runtimes of selected methods in seconds. Every box represents $45$ computed instances. \emph{Left:} $32 \times 32$. \emph{Right}: $64 \times 64$.}
  \label{fig:boxplots}    
\end{figure}

The runtimes of the tests are listed in Table~\ref{32table} for $32 \times 32$ and Table~\ref{64table} for $64 \times 64$, respectively, averaged over all $45$ instances in one class. The average over all classes can be found in the bottom row under 'Overall'. The fastest algorithm for each class is highlighted in bold. Additionally, boxplots for four selected methods are given in Figure~\ref{fig:boxplots}.

As the numbers show, the shielding neighborhood method is clearly the fastest algorithm for $32 \times 32$ instances among the methods tested. It takes hardly more than half the time on average compared to the transportation simplex, the shortlist method and the AHA method. The default solvers of CPLEX and Gurobi perform particularly badly. It is remarkable, however, that the network simplex solver of CPLEX outperforms the default solvers and the Gurobi dual simplex by a huge margin. The performance is still somewhat worse than our implementations of the transportation simplex and the shortlist method.

\begin{table}
\centering
\begin{tabular}{r||c|c|c|c|c|c|c|c|c|}
\multirow{2}{*}{Instance Class} & \multirow{2}{*}{TPS} & \multirow{2}{*}{SHL} & \multicolumn{2}{c|}{Shielding} & \multirow{2}{*}{AHA} & \multicolumn{2}{c|}{CPLEX} & \multicolumn{2}{c|}{Gurobi} \\
& & & CPX & TPS & & Def & NWS & Def & DS \\
\hline \hline
WhiteNoise & 1.58 & 1.38 & \textbf{0.67} &1.3 & 3.28 & 29.8 & 5.76 & 8.1 & 50.5 \\
\hline
GRFrough & 2.16 & 1.98 & \textbf{1.08} & 2.3 & 3.19 & 43.0 & 5.94 & 9.0 & 50.7 \\
\hline
GRFmoderate & 3.45 & 3.80 & \textbf{1.86} & 6.5 & 3.17 & 77.3 & 6.14 & 21.1 & 49.8 \\
\hline
GRFsmooth & 4.23 & 5.46 & \textbf{2.66} & 10.0 & 4.39 & 101.9 & 6.15 & 36.0 & 50.4 \\
\hline
LogGRF & 5.18 & 6.40 & \textbf{3.00} & 10.8 & 6.80 & 119.9 & 6.23 & 49.5 & 51.7 \\
\hline
LogitGRF & 4.50 & 5.26 & \textbf{2.40} & 7.9 & 8.49 & 98.4 & 6.33 & 31.0 & 51.7 \\
\hline
CauchyDensity & 4.46 & 6.06 & \textbf{3.62} & 12.7 & 3.76 & 140.4 & 6.09 & 54.1 & 49.3 \\
\hline
Shapes & 1.06 & 1.07 & \textbf{0.92} & 8.9 & 1.27 & 8.9 & 1.22 & 5.2 & 9.1 \\
\hline
ClassicImages & 3.33 & 3.26 & \textbf{1.58} & 5.2 & 2.01 & 68.5 & 6.05 & 18.0 & 49.6 \\
\hline
Microscopy & 2.34 & 3.02 & \textbf{1.66} & 9.4 & 3.14 & 35.2 & 2.74 & 20.7 & 22.1 \\
\hline \hline
Overall & 3.23 & 3.77 & \textbf{1.94} & 7.5 & 3.95 & 72.3 & 5.26 & 25.3 & 43.5 \\
\end{tabular}
\caption{Average runtimes on $32 \times 32$ instances in seconds. The columns represent the methods tested: transportation simplex (TPS), shortlist method (SHL), shielding neighborhood method with CPLEX (CPX) and TPS as internal solvers, AHA method, CPLEX with default (Def) parameters and the network simplex solver (NWS), as well as Gurobi with default parameters and the dual simplex solver (DS).}
\label{32table}
\end{table}

\begin{table}
\centering
\begin{tabular}{r||c|c|c|c|c|c|c|c|c|}
\multirow{2}{*}{Instance Class} & \multirow{2}{*}{TPS} & \multirow{2}{*}{SHL} & \multicolumn{2}{c|}{Shielding} & \multirow{2}{*}{AHA} & \multicolumn{2}{c|}{CPLEX} & \multicolumn{2}{c|}{Gurobi} \\
& & & CPX & TPS & & Def & NWS & Def & DS \\
\hline \hline
WhiteNoise & 74 & 56 & \textbf{13} & 65 & 36 & 2057 & 174 & 311 & 1657 \\
\hline
GRFrough & 153 & 110 & 24 & 153 & \textbf{20} & 3216 & 174 & 473 & 1659 \\
\hline
GRFmoderate & 261 & 306 & 51 & 469 & \textbf{23} & 4592 & 190 & 1971 & 1634 \\
\hline
GRFsmooth & 306 & 468 & 80 & 778 & \textbf{57} & 5621 & 195 & 3723 & 1590 \\
\hline
LogGRF & 439 & 531 & 79 & 756 & \textbf{59} & 7156 & 198 & 4628 & 1687 \\
\hline
LogitGRF & 333 & 362 & \textbf{69} & 552 & 77 & 6024 & 198 & 2294 & 1637 \\
\hline
CauchyDensity & 245 & 397 & 97 & 968 & \textbf{30} & 6336 & 195 & 4461 & 1435 \\
\hline
Shapes & 73 & 73 & 25 & 893 & \textbf{12} & 885 & 40 & 302 & 219 \\
\hline
ClassicImages & 218 & 246 & 41 & 418 & \textbf{18} & 6298 & 189 & 1551 & 1546 \\
\hline
Microscopy & 28 & 37 & 26 & 824 & 24 & 352 & \textbf{18} & 179 & 114 \\
\hline \hline
Overall & 213 & 258 & 51 & 588 & \textbf{36} & 4254 & 157 & 1989 & 1318 \\
\end{tabular}
\caption{Average runtimes on $64 \times 64$ instances in seconds. See the caption above for details.}
\label{64table}
\end{table}

At resolution $64 \times 64$ we see a similar picture with some exceptions. The shielding method (with CPLEX as internal solver) is even further ahead of most other algorithms, but at the same time the AHA method, which seems to be scaling much better than the linear programming approaches, has gained even more and in fact shows the best times now for many of the classes. However, keep in mind that the results of this method are not exactly correct and the timing varies according to the stopping criterion one applies (see the next subsection).

The CPLEX network solver is with the exception of the classes WhiteNoise and GRFrough quite a bit faster at resolution $64 \times 64$ than our implementations of the transportation simplex and the shortlist methods. In the two classes Shapes and Microscopy, which have many zeros in the images, it is even competitive with the shielding and AHA methods. Overall, it shows the most consistent performance both across the various benchmark classes (if effective size of the problem is taken into account) and within each class; see the last row in Figure~\ref{fig:boxplots}.

In contrast, the other methods show a much stronger sensitivity with respect to the class considered. Especially for classes~4--7 and to some extent also for class~3, we see higher average computation times and in particular a much wider spread of times including a number of outliers in Figure~\ref{fig:boxplots}.

\subsection{Errors of the AHA Method}
\label{errors}
As described in Subsection~\ref{aha}, the AHA method does not solve the problems we consider here with full accuracy, but makes a blurring error by interpreting pixel values of the source measure $\mu$ as uniformly distributed over small squares and a precision error by tackling a continuous minimization problem which for numerical reasons cannot be solved exactly. 

In Table~\ref{tab:errortable} we report the precision error (PE) in terms of the mass in the probability measure $\mu$ that is wrongly allocated,
as well as the relative Wasserstein error (RWE) made with the AHA method, i.e.
\begin{equation*}
  \frac{W_2^{\mathrm{AHA}}(\mu,\tilde{\nu})-W_2^{\mathrm{TSP}}(\mu,\nu)}{W_2^{\mathrm{TSP}}(\mu,\nu)},
\end{equation*}
where AHA and TSP denote the methods used and $\tilde{\nu}$ denotes the second marginal of the transference plan returned by AHA.

We can see that the precision error is reasonably small. What is more, if we assume that we would have to reroute the wrongly allocated mass roughly by a distance that corresponds to the true Wasserstein distance between $\mu$ and $\nu$ (this is reasonable in the sense that it is roughly the same order of magnitude as relevant distances measured in the image), we can compare the precision errors to the relative Wasserstein errors and see that the former play only a minor role. Consequently, the RWE is mainly due to the blurring effect. This is corroborated further by the fact that the RWE for the $64 \times 64$ resolution is considerably smaller than for $32 \times 32$.

\begin{table}[ht]
\centering
\begin{tabular}{r||c|c|c|c|}
 
 Instance Class & PE 32 & PE 64 & RWE 32 & RWE 64 \\ 
  \hline \hline
WhiteNoise & 0.3098e-05 & 0.3400e-05 & 0.031224 & 0.027232 \\ \hline
GRFrough & 0.4261e-05 & 0.6327e-05 & 0.018957 & 0.008058 \\ \hline
GRFmoderate & 0.9840e-05 & 1.8471e-05 & 0.003965 & 0.001075 \\ \hline
GRFsmooth & 2.4910e-05 & 4.0698e-05 & 0.001303 & 0.000349 \\ \hline
LogGRF & 4.6378e-05 & 5.6221e-05 & 0.000556 & 0.000230 \\ \hline
LogitGRF & 1.6833e-05 & 3.7144e-05 & 0.002008 & 0.000682 \\ \hline
CauchyDensity & 3.9273e-05 & 8.2546e-05 & 0.000988 & 0.000289 \\ \hline
Shapes & 1.3435e-05 & 3.3610e-05 & 0.003192 & 0.000994 \\ \hline
ClassicImages & 0.7797e-05 & 1.8311e-05 & 0.005119 & 0.001512 \\ \hline
Microscopy & 2.7111e-05 & 6.1215e-05 & 0.001243 & 0.000412 \\ \hline \hline
Overall & 1.9294e-05 & 3.5794e-05 & 0.006855 & 0.004083 \\ 
\end{tabular}
\caption{Average precision error (PE) and relative Wasserstein error (RWE) over the ten classes. See text for details.} 
\label{tab:errortable}
\end{table}

\subsection{Iterations of the Shielding Method}
\label{shielditer}
The shielding method solves the optimal transport problem via a sequence of restricted instances. Here we have a look on how many iterations of these instances are necessary.

\begin{figure}[th]
  \centering
  \includegraphics[width=\textwidth]{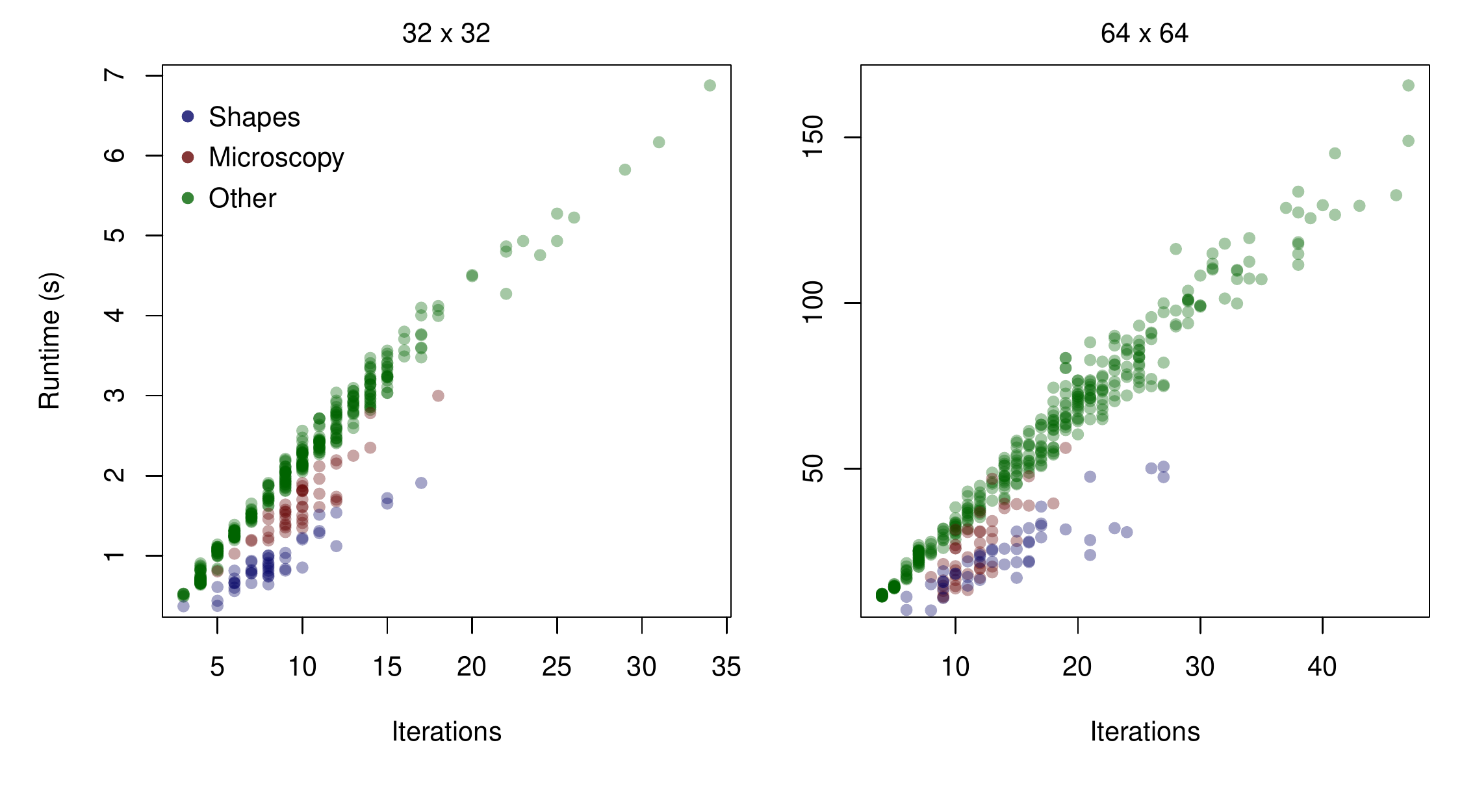}
  \vspace*{-11mm}

  \caption{Scatterplots for Shielding-CPX showing the runtimes against the number of iterations for the classes Shapes (blue), Microscopy (red) and the other classes combined (green). Every data point represents one of the $450$ computed instances. \emph{Left:} $32 \times 32$. \emph{Right:} $64 \times 64$.}
  \label{fig:shielditer}    
\end{figure}

On the scale $32 \times 32$ the average number of iterations varies between $3.9$ for WhiteNoise and $16.6$ for CauchyDensity. The numbers for scale $64 \times 64$ are higher ($4.6$ through $28.5$), but show similar behavior otherwise. Figure \ref{fig:shielditer} presents scatterplots of the runtime against the number of iterations.

For most classes (green points) we observe a linear scaling of the runtime with the number of iterations. This means that the runtime in each iteration is roughly the same over these classes and we may conclude, since CPLEX has runtimes that scale consistently with model size, that the neighborhood sizes remain more or less constant as well.
The only exceptions are the classes Shapes (blue) and Microscopy (red), where the runtimes are lower than expected from the number of iterations. This can be attributed to the internal solver, which benefits from the lower effective dimension that comes from the zero mass pixels occuring in these two classes. The difference is not as significant as for the global CPLEX network simplex runtimes, since the dimension has already been reduced by the construction of the shielding neighborhoods.

\section{Discussion}

If we look at the different classes, we note that the solvers perform much better on the classes where the number of pixels with mass zero is large (Microscopy and Shapes). This is because they seem to benefit particularly from the reduction of the effective dimension of the problem. The runtimes are very consistent across the other classes, which allows the conclusion that the solvers can only exploit the mathematical structure of the model, but are unable to use geometric features of the input data to their advantage.

Other linear programming methods benefit from the lower effective dimension as well, although the difference is not as significant. However, these methods seem to be comparatively faster on classes where most of the transports are rather short (rough structure, such as GRFrough or WhiteNoise), and slower on classes with longer transports (smooth structure, such as GRFsmooth or LogGRF). This can be explained by the initial solution routine, which is shared across many of the tested methods. The greedily selected initial transport plan, which favors short transports, is more likely to be close to optimal in short range classes than in long range classes. The shortlist method, which performs another greedy step in addition when searching for new basis variables within the shortlists, benefits particularly from short transports in the solutions, but not as much as in the sparse examples considered in \cite{GottschlichSchuhmacher2014} with the Euclidean metric as cost function.

The runtimes of the AHA method are relatively consistent. They are only considerably shorter for the class Shapes. This may be due to the fact that in these instances there are only a few different mass values in the images and the mass is uniformly distributed over large areas.

Interestingly, a comparison of the transportation simplex and the CPLEX network simplex reveals that the performance of the transportation simplex is better on $32 \times 32$ instances, while at resolution $64 \times 64$ the opposite is the case. This can be explained by looking at the two methods at hand. Although details of the CPLEX network simplex solver are not known, it is safe to assume that the simplex steps are implemented very efficiently. On the other hand, the transportation simplex has the advantage of an easily obtainable good initial solution, whereas in the network simplex a preceding simplex phase is necessary. This makes the TPS perform better on smaller instances. On the higher resolution our results suggest that the advantage of a strong initial solution is not as influential as the efficient simplex steps.

Considering the small disparity between the transportation simplex and the CPLEX network simplex runtimes, the difference in performance between the two internal solvers for the shielding neighborhood method is surprisingly large. This is due to the fact that initial solutions to the interior models are available in the shielding method and thus the first phase of the network simplex is not necessary. This is why the initial advantage of the transportation simplex disappears and hence using CPLEX as the internal solver yields much lower runtimes.

Another observation worth mentioning is that the runtimes, and therefore numbers of iterations, of the shielding method for the randomly created classes 1--6 agree very well with the ranges of dependence within the data of the classes. The class with the lowest runtimes, WhiteNoise, has no dependence at all, whereas the classes with long range dependences, GRFsmooth and LogGRF, belong to the more difficult classes for the shielding method to solve. That seems to indicate that constructing small shielding neighborhoods prevents larger updates to the current solution per iteration, which might be required in these instances. This may also contribute to the comparatively long and inconsistent runtimes on the class CauchyDensity.

Based on this first comparison of singlescale methods, we give the following \textbf{recommendations} for large discrete transport problems (on grids):
\begin{enumerate}
  \item If you have to be precise and have access to the IBM CPLEX software, use the shielding method in combination with CPLEX.\\[-2.5mm]
  \item If you have to be reasonably precise but can afford a small controllable error, use the AHA method. This is especially advisable if you require a very high resolution for the mass distribution at the source, as this comes for the AHA method at virtually no extra cost. \\[-2.5mm]
  \item Both methods are not widely available nor typically very efficient for costs other than the squared Euclidean metric. So for other costs direct use of a conventional simplex algorithm or the shortlist method may be preferable. \\[-2.5mm]
  \item If you use a solver directly, be very careful which one you use and that you call the most appropriate function. Especially when using CPLEX, make sure that you use the network simplex solver and that you set up the input as a network structure. If you can, solve the model with a warm start.
\end{enumerate}
\vspace*{1mm}

We emphasize again that the absolute runtimes given in Tables~\ref{32table} and \ref{64table} should not be taken at face value and that actual computations on modern CPUs are typically much faster. While the relative comparison presented here is justified to the best of our knowledge, it allows only limited conclusions about the performance of multiscale variants and multithreaded implementations of the different methods.

\section{Outlook}

By providing this benchmark we hope to improve the comparability of different methods for solving discrete optimal transport problems. Contributions or suggestions for extending the benchmark are welcome. In particular we plan to include data sets concentrated on more general grid structures and especially with irregular support points if there is enough public interest.

The {\sf R} package \texttt{transport} \cite{transport} offers user-friendly implementations that are mostly written in C of three of the methods presented here (transportation simplex, shortlist and the AHA method). It will be updated in the near future to include additional state-of-the-art methods.

Solving transport problems exactly for larger images (e.g.\ with $128$, $256$ or $512$ pixels in each dimension)
is still computationally very demanding, even for state-of-the-art methods. 
Efficient solutions of such large problems could pave the way for a new class of algorithms in image processing.
In the area of computer vision and image processing, important applications include 
image enhancement, denoising, inpainting, feature extraction and compression.
In one subdomain of image processing, these challenges are approached by decomposing images into two or three parts \cite{ThaiGottschlich2016DG3PD},
e.g. a cartoon component, which contains piecewise constant or piecewise smooth parts, 
a texture component, which captures oscillating patterns, 
and a noise component, which contains small scale objects (corresponding to high frequency parts in the Fourier domain).
After the decomposition step, the texture component can be utilized for applications such as fingerprint segmentation \cite{ThaiGottschlich2016G3PD}.
Image decompositions are obtained by formulating and solving minimization problems
that impose suitable norms on the respective components. 
The total variation (TV) norm is commonly used for the cartoon component and the G-norm \cite{Meyer2001}
for the texture component. 
Recent works by Brauer and Lorenz \cite{BrauerLorenz2015} and by Lellmann \textit{et al.} \cite{LellmannLorenzSchoenliebValkonen2014}
connect the G-norm to solutions of transport problems. 
Typically, the minimization problems described above are solved iteratively 
and are computationally expensive. 
It is conceivable to formulate transport norms for image decomposition,
which would require to solve a large transport problem in each iteration. 
Thus, efficient algorithms for optimal transport are a prerequisite for future research
in this direction.

\section*{Acknowledgements}
  The authors would like to thank Bj\"orn B\"ahre for implementing the AHA algorithm as part of a student project. They are also grateful to Stefan Stoldt and Stefan Jakobs (\url{https://jakobs.mpibpc.mpg.de/}) for sharing their microscopy images and for their permission to include these in the benchmark.

\bibliographystyle{spmpsci}

\bibliography{transport}

\end{document}